\documentclass[12pt]{article}
\begin{document}
\newcommand{\qed}{\ \
\mbox{\rule{8pt}{8pt}}\vspace{0.3cm}\newline}
\newcommand{\ia}{{\bf I}_{A_i}}
\newcommand{\ba}{\widetilde{{\bf KC}}_i}
\newcommand{\bba}{{\bf KC}_i}
\newcommand{\ra}{\longrightarrow}
\newcommand{\pe}{{\cal P}}
\newcommand{\pec}{\widehat{{\cal P}}}
\newcommand{\der}{{\cal DP}_d}
\newcommand{\adp}{{\cal DP}_d^{af}}
\newcommand{\ot}{\otimes}
\newcommand{\rec}{\raisebox{-1ex}{\ $\stackrel{\textstyle{\stackrel{\textstyle{\longleftarrow}}{\longrightarrow}}}{\longleftarrow}$\ }}
\newcommand{\Ain}{A_{\infty}}
\newcommand{\ka}{\mathbf{k}}
\newcommand{\bu}{\bullet}
\newcommand{\fo}{{\cal F}}
\newcommand{\pu}{{\cal P}}
\newcommand{\np}{{\bf N}[\frac{1}{p}]}
\title{On spectra and affine strict polynomial functors}
\author{Marcin Cha\l upnik\\
\thanks{The author was supported by the  Narodowe Centrum Nauki grants no. 2011/01/B/ST1/06184 and 2015/19/B/ST1/01150.}\\
%\thanks{The author was supported by the  Narodowe Centrum Nauki grant no. 2011/01/B/ST1/06184.}\\
\normalsize{Institute of Mathematics, University of Warsaw,}\\
\normalsize{ul.~Banacha 2, 02--097 Warsaw, Poland;}\\
\normalsize{e--mail: {\tt mchal@mimuw.edu.pl}}}
\date{\mbox{}}
%\thanks{The author was  supported by the grant  no. N N201 387034}}
\newtheorem{prop}{Proposition}[section]
\newtheorem{cor}[prop]{Corollary}
\newtheorem{theo}[prop]{Theorem}
\newtheorem{lem}[prop]{Lemma}
\newtheorem{defi}[prop]{Definition}
\newtheorem{defipro}[prop]{Definition/Proposition}
\newtheorem{fact}[prop]{Fact}
\newtheorem{exa}[prop]{Example}
\newcommand{\ca}{\mbox{$\cal{A}$}}
\newcommand{\la}{\lambda}
\maketitle
\begin{abstract}
We compare derived categories of the category of strict polynomial functors over a finite
field and the category of ordinary endofunctors on the category of vector spaces. We introduce two intermediate categories: the category of $\infty$--affine strict polynomial
functors and the category of spectra of strict polynomial functors. They provide a conceptual framework for computational theorems of Franjou--Friedlander--Scorichenko--Suslin and clarify the role of inverting the Frobenius
morphism in comparison between rational and discrete cohomology.
\\\mbox{}\vspace{0.1cm}
\end{abstract}
\section{Introduction}
The aim of the present paper is to better understand relationship between derived categories of the category ${\cal P}_d$ of strict polynomial functors of degree $d$ 
over a finite field $\ka$  and the category ${\cal F}$ of usual functors on
the vector spaces over $\ka$.
It may be thought of as an instance of a fundamental problem in algebraic geometry: comparing  affine schemes  with their sets of rational points over  small
fields. In our context, the situation
is well understood at the level of abelian categories. Namely, it was shown  [FFSS, Prop. 1.4] that if $d\leq |\ka|$, then the forgetful
functor ${\bf f}:{\cal P}_d\ra{\cal F}$ is a full embedding.  Unfortunately, this is no longer true at he level of derived categories, as
the forgetful functors does not induce an isomorphism on Ext--groups. However, as it was shown by Franjou--Friedlander--Scorichenko--Suslin, we still get an isomorphism
on Ext--groups when we take instead of given strict polynomial functors $F,G$ their
large enough Frobenius twists $F^{(ni)},G^{(ni)}$. To put it more precisely: since
${\bf f}(F)={\bf f}(F^{(ni)}), {\bf f}(G)={\bf f}(G^{(ni)})$ for $|\ka|=p^n$, we have the induced map
\[{\rm colim}_i\ {\rm Ext}^*_{{\cal P}_{dp^{ni}}}(F^{(ni)},G^{(ni)})\ra {\rm Ext}^*_{{\cal F}}(F,G),\]
and [FFSS, Th. 3.10] says that this map is an isomorphism provided that
 $d\leq |\ka|$.
 This result has numerous nontrivial applications, since in many cases the left hand side is much more computable. Unfortunately, neither its applications nor
 its proof answers a natural question: why does twisting make
 ${\cal P}_d$ and ${\cal F}$ closer to each other? \newline In our article we  address this question by putting 
 the above mentioned phenomena into a wider categorical context (although it should be emphasized that our work does not provide provide new proofs of the results of [FFSS] , since we use their computations in several places). Namely, we factorize
 the (derived)  forgetful functor ${\bf f}: {\cal DP}_d\longrightarrow {\cal DF}$ through certain intermediate triangulated category whose Hom-spaces are (among others) colimits of Ext--groups between twists of strict polynomial functors. In fact, we construct two apparently quite different triangulated 
 categories which have this property. \newline The first construction which is described in Sections 2 and 3 uses a concept of ``affine strict polynomial functor'' introduced in [C4]. In fact our situation resembles that considered in [C4]. 
 In both cases we would like to extend a reflective full embedding of abelian categories to their derived categories but we face an obstruction that the unit of our adjunction is not an isomorphism. 
 Nevertheless, this unit admits an explicit description which gives us a hint how to enlarge 
 the starting category to obtain a reflective full embedding (in fact, as shown in [Ku2] the full embedding  ${\bf f}$ is a part of recollement diagram of abelian categories, however in our situation, 
 similarly to [C4], we do not get
 a full recollement due to appearance of categories of infinite homological dimension).
 By applying this procedure we get a DG functor category ${\cal P}_d^{af_{\infty}}$ and then
 its derived category ${\cal DP}_d^{af_{\infty}}$ together with a factorization of ${\bf f}$ as
 \[{\cal DP}_d\stackrel{z^*}{\ra}{\cal DP}_d^{af_{\infty}}\stackrel{{\bf f}^{af_{\infty}}}{\ra}
 	{\cal DF}.\]
The  important features of this factorization are that 
 \[{\rm Hom}^*_{{\cal DP}_d^{af_{\infty}}}(z^*(F),z^*(G))\simeq {\rm colim}_i
 {\rm Ext}^*_{{\cal P}_{dp^{i}}}(F^{(i)},G^{(i)})\]
 and that ${\bf f}^{af_{\infty}}$ (when restricted to the subcategory ${\cal DP}_d^{faf_{\infty}}$ 
of finite objects in ${\cal DP}_d^{af_{\infty}}$) is a full embedding (Theorem 3.7).  Thus we have achieved
our goal by a rather tautological construction, since the category ${\cal DP}_d^{af_{\infty}}$ is designed exactly in such a way that we obtain the desired colimits as Ext-spaces. This being
said, the fact that the construction works is highly nontrivial because it utilizes in an essential way the
formality phenomena observed in [C3, C4]. \newline  Then in Sections 4 and 5 we take
quite a different approach, which is perhaps more intuitive. It relies on an observation that an important difference between the categories ${\cal P}$ and ${\cal F}$ is that in the latter the
Frobenius twist operation is invertible. Hence if we formally invert the Frobenius twist in ${\cal P}$ we should obtain a category closer to ${\cal F}$. Moreover, when we think of classical example of applying such a construction i.e. stable homotopy category, we see how 
colimits enter our story: we should get them as an analog of the known description of the homotopy classes of maps between suspension spectra.\newline Technically, we introduce the category ${\cal SP}_d$ of spectra 
of complexes of strict polynomial functors and, following a general approach of Hovey [Ho1], we introduce a  Quillen model structure on it. Then we put  ${\cal DSP}_d$
to be the homotopy category with respect to this structure and we get a factorization of 
${\bf f}$ as
\[{\cal DP}_d\stackrel{{\bf C}^{\infty}}{\ra}{\cal DSP}_d\stackrel{{\bf f}^{st}}{\ra}
{\cal DF},\]
where ${\bf C}^{\infty}$ is a functor analogous to $\Sigma ^{\infty}$ in topology.
Then we have the expected description of the maps between ``suspension spectra'': (Theorem 4.6): 
\[{\rm Hom}^*_{{\cal DSP}_d}({\bf C}^{\infty}(F),{\bf C}^{\infty}(G))\simeq {\rm colim}_i
{\rm Ext}^*_{{\cal P}_{dp^{i}}}(F^{(i)},G^{(i)})\]
 and, similarly to the first approach, ${\bf f}^{st}$ when restricted to the category ${\cal 
 	DP}_d^{st}$ generated as triangulated category with direct summands by the image  of ${\bf C}^{\infty}$, is a full embedding (Theorem 5.2).\newline
 Finally, in the last section we compare the both constructions. Namely, we find a full
 embedding  $\gamma:
 {\cal DP}^{af_{\infty}}_d\ra {\cal DSP}_d$ which restricts to an equivalence 
 ${\cal DP}_d^{faf_{\infty}}\simeq {\cal DP}_d^{st}$ (Theorem 6.2). This shows that the categories 
${\cal DP}^{af_{\infty}}_d$ and ${\cal DSP}_d$ are quite close, which is
perhaps not obvious  at a first sight. As a  sort of heuristic explanation we can offer the following observation.
Similarly to  the classical context, in  our category of spectra, the delooping functor $\Theta^{\infty}$ plays an important role. On the other hand, as it was observed in [C4],
the category of affine strict polynomial functors is closely related to the category of representations of the group of algebraic loops on $GL_n(\ka)$. Thus our category 
${\cal P}_d^{af_{\infty}}$ should correspond to the infinite loops on $GL_n(\ka)$. Hence
some sort of relation to infinite loop spaces is a feature shared by the both 
categories.\newline
Now let us discuss the differences between ${\cal DP}_d^{af_{\infty}}$ and ${\cal DSP}_d$.
The fact that ${\cal DP}_d^{af_{\infty}}$ embeds into ${\cal DSP}_d$ shows that
the former category is closer to ${\cal DP}_d$. This is not surprising, since we see in its
very construction, that it is a possibly closest to ${\cal DP}_d$ triangulated category in which the colimits
of Exts of twists appear (we make no attempt to make this statement precise). In particular, we see that it is not necessary to fully invert the Frobenius twist to
get these colimits. In fact, one can show that the Frobenius twist gives a full embedding
${\cal DP}_d^{af_{\infty}}\subset {\cal DP}_{dp}^{af_{\infty}}$ but not an equivalence.\newline
Thus, one could think that the factorization through ${\cal DSP}_d$ is something less
fundamental. On the other hand however, the functor ${\bf f}^{st}: {\cal DSP}_d\ra {\cal DF}$ has a remarkable property that it preserves (at least some) fibrant objects (Remark 5.3).
This suggests that ${\cal SP}_d$ (in contrast to just ${\cal P}_d$) encodes important information about injective objects in
${\cal F}$.  \newline
This article was inspired by an ongoing project joint with Piotr Kowalski exploring
connections between rational cohomology and difference algebra. I am grateful to
Piotrek for many discussions on various aspects of Frobenius twist. I would also 
like to thank Stanis\l aw Betley for remarks on the preliminary version of this article.\newline
In order to help the reader navigating in the article I provide below the list of main definitions and notations used in the paper.

\vspace{1cm}

\begin{tabular}{l|c|l}
notation & Section & meaning\\
\hline
$\ka$ & 2 & ground field, since Section 3 finite\\
$\Ain$ & 2 & graded algebra $\ka[x_1,x_2,\ldots]/(x_1^p,x_2^p,\ldots)$\\
${\cal V}$ & 2  & category of finite-dimensional vector spaces over $\ka$\\
${\cal V}'$ & 2  & category of all vector spaces over $\ka$\\
${\cal V}^{gr}$ & 2  & category of graded vector spaces over $\ka$ finite-dimensional\\
&& in each degree\\
${\cal V}^{grf}$ & 2  & category of totally finite-dimensional graded vector spaces over $\ka$ \\
${\cal V}'^{gr}$ & 2  & category of all graded vector spaces over $\ka$\\
${\cal V}_{\Ain}$ & 2& subcategory of the category of free graded $\Ain$-modules\\ 
$\Gamma^d{\cal V}^{gr}$ & 2 & category of divided powers over ${\cal V}^{gr}$\\
$\Gamma^d{\cal V}^{grf}$ & 2 & category of divided powers over ${\cal V}^{grf}$\\
$\Gamma^d{\cal V}_{\Ain}$ & 2 & category of divided powers over ${\cal V}_{\Ain}$\\
${\cal P}_d$ & 2 & category of strict polynomial functors over $\ka$ of degree $d$\\
${\cal P}_d^{gr}$ & 2 & category of graded strict polynomial functors over $\ka$ of degree $d$\\
${\cal P}_d^{af_{\infty}}$ & 2 & category of $\infty$-affine strict polynomial functors of degree $d$\\
$z^*$ & 2 & functor from ${\cal DP}_d$ to ${\cal DP}_d^{af_{\infty}}$ induced by the forgetting\\
 $t^*$ & 2 & functor from ${\cal DP}_d^{af_{\infty}}$ to ${\cal DP}_d$ induced by the scalar extension,\\
 && right adjoint to $z^*$\\
$\Gamma^{d,U}$ & 2 &object in ${\cal P}_d$ represented by $U\in\Gamma^d{\cal V}$\\
$h^{U\ot\Ain}$ & 2 &object in ${\cal P}_d^{af_{\infty}}$ represented by 
$U\ot\Ain\in\Gamma^d{\cal V}_{\Ain}$\\
${\cal F}$ & 3 & category of functors from ${\cal V}$ to ${\cal V}'$\\
${\bf f}$ &3&  forgetful functor from ${\cal DP}_d$ to ${\cal DF}$ \\
%${\bf Ra}$ &3& functor from ${\cal DF}$ to ${\cal DP}_d$ rght adjoint to ${\bf f}$\\
${\cal DP}_d^{faf_{\infty}}$ & 3 & smallest triangulated subcategory of 
${\cal DP}_d^{af_{\infty}}$containing $h^{U\ot\Ain}$\\
&& and closed under isomorphisms and direct summands\\
${\bf C}$ &4& Frobenius twist regarded as functor from 
${\cal DP}_d$ to ${\cal DP}_{pd}$\\
$\widehat{{\cal P}}$ & 4 & product category $\prod_{d>0}{\cal P}_d$\\
${\cal K}{\cal P}_d$ & 4& category of complexes over ${\cal P}_d$\\
${\cal K}\widehat{{\cal P}}$ & 4& category of complexes over $\widehat{{\cal P}}$\\
${\cal S}\widehat{{\cal P}}$ & 4& category of spectra over ${\cal K}\widehat{{\cal P}}$\\
${\cal S}{\cal P}_d$ & 4& category of spectra over ${\cal KP}_d$\\
${\cal K}{\cal F}$ & 5& category of complexes over ${\cal F}$\\
${\cal S}{\cal F}$ & 5& category of spectra over ${\cal KF}$\\
\end{tabular}
\newpage
\begin{tabular}{l|c|l}
${\cal DSP}_d^{st}$ &5 & smallest triangulated subcategory of ${\cal DSP}_d$ containing
${\bf C}^{\infty}F)$\\
&& and closed under isomorphisms and direct summands\\
${\cal P}_d^{af_i}$ & 6 & category of $i$-affine strict polynomial functors of degree $d$\\
\end{tabular}

\section{$\infty$--affine functors}
In this section we introduce and establish basic properties of the category 
${\cal P}_d^{af_{\infty}}$ of $\infty$--affine strict polynomial functors. 
In the next section we relate ${\cal P}_d^{af_{\infty}}$ to the categories 
${\cal P}_d$ and ${\cal F}$.\newline 
We shall frequently use graded categories and functors. There are several possible choices 
of the setup here, we follow the one from [Ke] and we briefly recall the basic definitions. We fix a ground  field $\ka$. By a graded ($\ka$-linear) category we mean a $\ka$-linear category with Hom-spaces equipped with ${\bf Z}$-grading preserved by composition, which means that if $|f|=n$, $|g|=m$ then $|g\circ f|=n+m$. 
The basic  example of graded $\ka$-linear category is the category 
${\cal V}'^{gr}$  of ${\bf Z}$--graded vector spaces over $\ka$. We put decoration $(-)'$ because in the majority of our constructions we will restrict to its subcategory ${\cal V}^{gr}$
consisting of graded vector spaces finite dimensional in each degree.
 The Hom-spaces are given as:
 \[\mbox{Hom}^n_{{\cal V}'^{gr}}(V^{\bullet},W^{\bullet}):=
\prod_j \mbox{Hom}(V^{j},W^{j+n})\]
% with the grading given by the following formula. For $f: V^i\longrightarrow W^j$ we 
 %have $|f|=j-i$.  
A graded functor between graded categories is a functor which preserves the grading on Hom-spaces.  
Observe that for any graded $\ka$-linear category ${\cal A}$ the  graded functors
from ${\cal A}$ to ${\cal V}'^{gr}$ form a $\ka$-linear graded category. In order to define the graded Hom-spaces let us define for a graded functor $F$ from ${\cal A}$ to ${\cal V}'^{gr}$
its shift $F[1]$ by $F[1](A)^n:=F(A)^{n+1}$. Then we define ${\rm Hom}^n(F,G)$ to be
${\rm Nat}(F,G[n])$ (see [Ke, Sect. 1]).\newline
The notion of $\infty$--affine strict polynomial functor is quite straightforward generalization of that of the affine strict polynomial functor introduced in [C4]. The only difference is that instead the graded algebra $A=\ka[x]/x^p\simeq 
\mbox{Ext}^*_{{\cal P}_{dp}}(I^{(1)},I^{(1)})$ we consider the graded algebra
$\Ain:=\ka[x_1, x_2,\ldots]/(x_1^p, x_2^p,...)$ with $|x_i|=2p^i$.  
The appearance of the algebra $\Ain$ in our situation follows from the fact known  from [FLS] that if   $\ka$ is a finite field of characteristic $p$ then there is   an isomorphism of graded algebras 
$\Ain\simeq \mbox{Ext}^*_{{\cal F}}(I,I)$. \newline
Let ${\cal V}$ (resp. ${\cal V}'$) stands for the category of  finite dimensional vector spaces over $\ka$
(resp. the category of all vector spaces over $\ka$).  Next, let ${\cal V}_{\Ain}$ stand for 
the full subcategory of the category of graded free $\Ain$-modules consisting of 
the modules $V\ot\Ain$ for $V\in{\cal V}$.
Then we consider the category
  $\Gamma^d{\cal V}_{\Ain}$ of divided powers over ${\cal V}_{\Ain}$ (see e.g. [FP, Section 3]).
  This is  the graded $\ka$-linear category with the objects the same as those
of ${\cal V}_{\Ain}$, but the  Hom-spaces are
\[\mbox{Hom}_{\Gamma^d {\cal V}_{\Ain}}(V\ot \Ain,W\ot \Ain):=\Gamma^d(\mbox{Hom}_{\Ain}(V\ot \Ain,W\ot\Ain))\]
where $\Gamma^d$ stands for the space of symmetric $d$-tensors over $\ka$. The grading comes from that  on $\Ain$ and the standard grading on the tensor product:
\[|x_1\ot\ldots\ot x_d|:=\sum_{s=1}^d |x_s|.\]
%We observe that the latter space can also be interpreted as the space of polynomial (over $%\ka$) maps homogeneous of degree $d$ from 
The composition law is given by the composition of the map:
\[
\Gamma^d(\mbox{Hom}_{\Ain}(V,W))\ot \Gamma^d(\mbox{Hom}_{\Ain}(W,U))\longrightarrow
\]
\[
\Gamma^d(\mbox{Hom}_{\Ain}(V\ot \Ain,W\ot \Ain)\ot 
\mbox{Hom}_{\Ain}(W\ot \Ain,U\ot \Ain))\]
coming from the natural map $\Gamma^d(X)\ot\Gamma^d(Y)\longrightarrow
\Gamma^d(X\ot Y)$ existing for any vector spaces $X,Y$, with the map
\[\Gamma^d(\mbox{Hom}_{\Ain}(V\ot \Ain,W\ot \Ain)\ot 
\mbox{Hom}_{\Ain}(W\ot \Ain,U\ot \Ain))\longrightarrow
\]
\[
\Gamma^d(\mbox{Hom}_{\Ain}(V\ot \Ain,U\ot \Ain))\]
coming from the composition in $\Ain$-linear Homs. 
 \begin{defipro}
 	An $\infty$--affine strict polynomial functor $F$ homogeneous of degree $d$ is a graded $\ka$--linear functor
 	\[F:\Gamma^d {\cal V}_{\Ain}\ra {\cal V}'^{gr}.\]
 	The affine strict polynomial functors homogeneous of degree $d$ with morphisms being natural transformations (with shifted targets) form a graded abelian category ${\cal P}_d^{af_{\infty}}$ (see [C4, pp. 655-656]).
 \end{defipro}
The reader of [C4] will find there a similar construction. In fact our category 
${\cal P}_d^{af_{\infty}}$ may be thought of as obtained by infinitely times applying the construction producing ${\cal P}_d^{af}$ from [C4]. 
However,   
we alert the reader that, in contrast to [C4] and also to the foundational paper 
on strict polynomial functors [FS], we do not impose any finiteness/finite generation assumptions
on values of functors.\newline 
Now we list some basic properties of the category ${\cal P}_d^{af_{\infty}}$.
The proofs of respective facts on affine functors from [C4] in most cases carry over to the current context. Therefore  we discuss more thoroughly  these points only when 
the infinite dimension of $\Ain$ requires special attention. Also, the reader interested in obtaining  more motivation behind the construction of affine functors is referred to [C4].\newline
Like in any functor category, for any $U\in{\cal V}$ we have the representable functor $h^{U\ot \Ain}\in {\cal P}_d^{af_{\infty}}$ given by the formula
\[V\ot \Ain\mapsto \mbox{Hom}_{\Gamma^d {\cal V}_{\Ain}}(U\ot \Ain,V\ot \Ain)\]
and the co--representable functor $c_{U\ot \Ain}^*\in {\cal P}_d^{af_{\infty}}$ given by the formula
\[V\ot \Ain\mapsto \mbox{Hom}_{\Gamma^d {\cal V}_{\Ain}}(V\ot \Ain,U\ot \Ain)^*\]
%\[V\ot \Ain\mapsto \mbox{Hom}_{\Gamma^d {\cal V}_{\Ain}}(V\ot \Ain,U\ot \Ain)^*\]
where $(-)^*$ stands for the graded $\ka$--linear dual. We list the basic properties of 
${\cal P}_d^{af_{\infty}}$. In part 3 when talking about projectives/injectives we regard ${\cal P}_d^{af_{\infty}}$ as just an abelian category (we forget an extra structure of grading).

\begin{prop}
	\begin{enumerate}
		\item
		There are  natural in $U\ot \Ain$ isomorphisms
		\[\mbox{Hom}_{{\cal P}_d^{af_{\infty}}}(h^{U\ot \Ain},F)\simeq F(U\ot \Ain)\]
		\[\mbox{Hom}_{{\cal P}_d^{af_{\infty}}}(F,c_{U\ot \Ain}^*)\simeq (F(U\ot \Ain))^*\]
		for any $F\in{\cal P}_d^{af_{\infty}}$.\newline
		\item
		Moreover, the map $\Psi: h^{U\ot \Ain}\otimes F(U\ot\Ain)\ra F$ adjoint to the map $F_{U\ot \Ain,V\ot \Ain}$ giving the action of $F$ on morphisms is surjective, provided that $\dim(U)\geq d$.
		\item
		If $\dim(U)\geq d$ then  $h^{U\ot \Ain}$ is a projective generator of ${\cal P}_d^{af_{\infty}}$, $c_{U\ot \Ain}^*$ is an injective generator of ${\cal P}_d^{af_{\infty}}$
\end{enumerate}
\end{prop}
The proofs of [C4, Prop. 2.5, 2.6] carry over to the current situation.\newline
Now we turn to comparing ${\cal P}_d^{af_{\infty}}$ with ${\cal P}_d$  (or rather its  graded variant ${\cal P}_d^{gr}$ which will be defined later) which, when compared to [C4], is a bit more delicate point due to infinite dimension of $\Ain$.
We have a pair of adjoint functors between the source categories of our functor categories. The first is just the forgetful functor:
\[z:\Gamma^d{\cal V}_{\Ain}\ra \Gamma^d{\cal V}^{gr}\]
($z$ comes from the polish word for forgetting which is {\em zapominanie}), where
$\Gamma^d{\cal V}^{gr}$ stands for the graded category whose objects are graded vector spaces 
finite dimensional in  each degree and 
\[\mbox{Hom}_{\Gamma^d{\cal V}^{gr}}^{\bullet}(V^{\bullet},W^{\bullet}):=
\Gamma^d(\mbox{Hom}(V^{\bullet},W^{\bullet}))\].
The second is the scalar extension functor:
\[t:\Gamma^d{\cal V}^{gr}\ra \Gamma^d{\cal V}_{\Ain}\]
explicitly given on objects as tensoring over $\ka$ with $\Ain$: $t(V):=V\ot \Ain$.
%$\Gamma^d{\cal V}^{gr}$ here is a variant of the category $\Gamma^d{\cal V}$
%where the underlying vector category is the category ${\cal V}^{gr}$ of ${\bf %Z}$--graded
%vector spaces.  
Then we claim that precomposing with 
$z$ can be extended to a   functor \[z^*:{\cal P}_d^{gr}\ra {\cal P}_d^{af_{\infty}}.\] 
where ${\cal P}_d^{gr}$
is the category of ${\bf Z}$--graded strict polynomial functors of degree $d$ (i.e. the category of graded functors from $\Gamma^d{\cal V}^{gr}$ to ${\cal V}'^{gr}$).. 
This can be done in two steps. The first step was already used in [C4, Sec. 2] (see also [T2, Sec. 2.5]).
Let $\Gamma^d{\cal V}^{grf}$ stands for the 
full subcategory of  $\Gamma^d{\cal V}^{gr}$ consisting of totally finite dimensional vector spaces.
Then any
$F\in{\cal P}_d$ can be extended to a graded
functor $F^{grf}:\Gamma^d{\cal V}^{grf}\ra {\cal V}'^{gr}$ (see [C4, p. 657]).
In the second step we shall extend  $F^{grf}$  to the graded functor $F^{gr}$ defined on the whole graded category  $\Gamma^d{\cal V}^{gr}$.  
For $V^{\bullet}\in{\cal V}^{gr}$ let $V^{\leq |j|}:=
\bigoplus_{|s|\leq j} V^s$. Then we define $F^{gr}(V^{\bullet})$ as $\mbox{colim}_j F^{grf}(V^{\leq |j|})$.
Now we can correctly define $z^*: {\cal P}_d^{gr}\ra {\cal P}_d^{af_{\infty}}$
%where ${\cal P}_d^{gr}$ means the graded category of graded functors from 
%$\Gamma^d{\cal V}^{gr}$ to ${\cal V}'$
by putting
\[z^*(F)(V\ot\Ain):=F^{gr}(V\ot \Ain).\]
Analogously, precomposing with $t$ produces the functor
\[t^*:{\cal P}_d^{af_{\infty}}\ra {\cal P}_d^{gr}.\]
Considering here the category ${\cal P}_d^{gr}$ instead of ${\cal P}_d$ will
be essential in the next section where we will compare the derived categories
of ${\cal P}_d$ and ${\cal P}_d^{af_{\infty}}$.\newline 
We list the properties of $z^*$ and $t^*$:
\begin{prop}\mbox{}\\
\begin{enumerate}
\item
	$z^*$ preserves representable objects i.e. \[z^*(\Gamma^{d,U})=h^{U\ot\Ain}.\]
	where $\Gamma^{d,U}\in{\cal P}_d$ is defined as $V\mapsto\mbox{Hom}_{\Gamma^d{\cal V}}(U,V)=\Gamma^d(U^*\ot V)$.
	\item
	The functor $t^*$ is right adjoint to $z^*$.
	\end{enumerate}
\end{prop}
Again, the proof of [C4, Prop. 2.4]
carries over to the present situation.\\
{\bf Remark: } It is worth mentioning that, as pointed out by the referee, all results of the present section hold for any non-negatively graded algebra finite dimensional in each degree put instead of $A_{\infty}$. In fact we do consider in Section 6 yet  another (and much simpler) variant of this construction. What is  specific to $A_{\infty}$ is relation of ${\cal P}_d^{af_{\infty}}$ to ${\cal P}_{pd}$ at the level of derived categories. We will discuss this relation in the next section.
\section{Formality and  $\infty$--affine algebraification} 
In the present section we factorize the derived functor of the forgetful functor ${\bf f}: {\cal P}_d\ra
{\cal F}$ through the derived category of  ${\cal P}_d^{af_{\infty}}$. Starting from this section  we assume that the ground field $\ka$ is a finite field of characteristic $p$.
As we have mentioned in Introduction we closely follow the strategy taken in [C4].\newline
To this end we  regard $\Gamma^d{\cal V}_{\Ain}$
as DG category with trivial differentials. Our main reference for general facts and terminology on DG categories is still [Ke].
We consider the category
$\mbox{Dif}(\Gamma^d{\cal V}_{\Ain}^{op})$ consisting of DG functors from
$\Gamma^d{\cal V}_{\Ain}$ to the category of complexes of $\ka$--modules (our strange terminology here is coherent with [Ke, Sect. 1.2], where the main focus was on contravariant functors). We are ready to introduce our main object of interest in the first part of the paper.
\begin{defi}
 Let
${\cal DP}_d^{af_{\infty}}$ be the category obtained from $\mbox{Dif}(\Gamma^d{\cal V}_{\Ain}^{op})$ by localization with respect to the class of quasi-isomorphisms
(we recall again that  we do not make any boundedness/finiteness assumptions). 
\end{defi}
 Let $\Gamma^d(I^*\ot I)$
denote  bifunctor given by the formula $(V,W)\mapsto \Gamma^d(V^*\ot W)$.
We regard $\Gamma^d(I^*\ot I)$ as a contravariant strict polynomial functor of degree $d$  in $V$ and just a naive functor in $W$. 
We shall denote the category of such mixed bifunctors by ${\cal P}_{\cal F}^d$ (we put decoration $d$ as a superscript to emphasize that our functors are contravariant with respect to the strict polynomial variable).
%Then by the Yoneda lemma the assignment
%\[F\mapsto \mbox{Hom}_{{\cal P}_d}(\Gamma^d(I\ot I^*),F)\]
%is nothing but the forgetful functor ${\bf f}:{\cal P}_d\ra {\cal F}$,
The assignment
\[F\mapsto \mbox{Hom}_{{\cal F}}(\Gamma^d(I^*\ot I),F)\]
defines the functor ${\bf a}:{\cal F}\ra {\cal P}_d$ which we call the (right) algebraification.
Then it is easy to see that ${\bf a}$ is right adjoint to the forgetful functor ${\bf f}:{\cal P}_d\ra {\cal F}$ (in fact, this adjunction is among
those considered by Kuhn in [Ku2]). %who has shown that it forms a part of %recollement diagram
%of abelian categories. 
Next, since ${\bf f}$ is an exact functor between abelian categories, it extends degreewise to a functor between unbounded derived categories
\[{\bf f}:{\cal DP}_d\ra {\cal DF}\]
which we shall denote by the same letter, though, formally it is the derived functor of ${\bf f}$. We shall abuse notation in such a manner in several further places in the article, where functors between abelian categories 
thanks to exactness simply factorize to derived category.
On the other hand ${\bf a}$ is a left exact  functor between abelian categories, hence it gives rise to a derived functor
\[{\bf Ra}:{\cal DF}\ra {\cal DP}_d.\]
Then it is a matter of routine verification that ${\bf f}$ and ${\bf Ra}$ remain adjoint (compare [C3, Th. 2.2]).
Our goal is to factorize this adjunction through ${\cal DP}_d^{af_{\infty}} $.\newline 
A crucial tool for studying the formality phenomena in our context is the following DG counterpart of the graded category $\Gamma^d{\cal V}_{\Ain}$.
\begin{defi}
Let $X^{\bullet}$ be a projective resolution of $\Gamma^d(I^*\ot I)$ in the category ${\cal P}^d_{\cal F}$.
We introduce the DG category $\Gamma^d{\cal V}_X^{op}$ whose objects are finite $\ka$-vector spaces and
\[{\rm Hom}^n_{\Gamma^d{\cal V}_X^{op}}(V,V'):={\rm Hom}_{{\cal F}}(X^{\bullet}(V',-),X^{\bullet}(V,-)[n]),\]
where $\rm{Hom}_{{\cal F}}$ stands for the Hom complex (ie. we do not require that the maps preserve differentials).
The composition law comes from the composition in Hom complexes.
\end{defi}
The connection between ${\cal P}_d^{af}$ and  $\Gamma^d{\cal V}_X^{op}$ over a field large enough essentially follows from the classical Ext-computations of Franjou-Friedlander-Scorichenco-Suslin:
\begin{prop}
Assume that $|\ka|\geq d$. Then the assignment $V\ot\Ain\mapsto V$  extends to a quasi-isomorphism  of  graded  categories $\Gamma^d{\cal V}_{\Ain}\simeq H^*(\Gamma^d{\cal V}_X^{op})$.
\end{prop}
{\bf Proof:} We need to establish a natural in $V,V'$ isomorphism of graded spaces:
\[{\rm Ext}^*_{{\cal F}}(\Gamma^{d,V'}, \Gamma^{d,V})\simeq \Gamma^d({\rm Hom}(V,V')\otimes\Ain).\]
In fact, for $V=V'=\ka$ it is just [FFSS, Th.6.3.(3)]. The general case may be obtained by a similar reasoning. Alternatively,  we can also  quickly derive Proposition 3.3 from ``the derived Kan extension''. Namely,
from [C3, Cor. 3.7] and the Yoneda lemma we get:
\[{\rm Ext}^*_{{\cal P}^{dp^i}}(\Gamma^{d(i),V'}, \Gamma^{d(i),V})\simeq 
\Gamma^d({\rm Hom}(V,V')\otimes\Ain).\]
Then our assertion follows from [FFSS, Th. 3.10]. \qed 
 However,  the much stronger fact, 
which generalizes [C4, Th. 4.2], holds:
\begin{theo}
	Assume that the ground field $\ka$ has $q\geq d$ elements.
	Then the assignment $V\ot\Ain\mapsto V$  extends to a quasi-isomorphism  of  DG categories $\phi: \Gamma^d{\cal V}_{\Ain}\simeq \Gamma^d{\cal V}_X^{op}$.
\end{theo}
The proof is conceptually  similar to that of [C4, Th. 4.2], but since there are some  technical differences we present it in some detail. 
First of all we need a certain generalization of Touz\'e universal classes [T1, T3].
\begin{lem}
There  exist  classes 
\[c[d]^{(i)}\in \mbox{Ext}^{2dp^{i-1}}_{{\cal P}_{dp^i}^{dp^i}}(\Gamma^{dp^i}(I^*\ot I),\Gamma^d(I^{*(i)}\ot
	I^{(i)}))\]
such that $c[1]^{(i)}\neq 0$ for all $i\geq 1$, and are compatible with cup product i.e.
\[\Delta_*(c[d]^{(i)})=(c[1]^{(i)})^{\cup d}\]
where $\Delta:\Gamma^d\ra I^d$ is the standard embedding.
\end{lem}  
{\bf Proof of Lemma 3.5} For $i=1$ we have the original Touz\'e classes, but the proof
carries over to this more general case. Indeed: it immediately follows from the degeneracy 
of the twisting spectral sequence  [T3, Prop. 17] and this degeneracy was 
showed also for multiple twists [T3, Th. 4]. \qed
Let $\Gamma^d( I^*\otimes I))_{\Ain})$ stands for the graded bifunctor $(V,W)\mapsto
\Gamma^d( V^*\otimes W\otimes\Ain)$ with grading coming from that on $A_{\infty}$ and 
$\Gamma^d$, regarded as an object in the derived category
of ${\cal P}^d_{{\cal F}}$.
Now we get an ${\cal F}$-analog of [C3, Prop. 3.3]. 
\begin{lem}
	There exist classes \[\widetilde{e}_d\in \mathrm{Hom}_{{\cal DP}^d_{{\cal F}}}(\Gamma^d(I^*\otimes I),(\Gamma^d( I^*\otimes I))_{\Ain^*})\] satisfying:
	\begin{enumerate}
		\item $\widetilde{e}_1\in\mathrm{Hom}_{{\cal DP}^1_{{\cal F}}}(I^*\otimes I,I^*\ot I\ot\Ain^*)$ is nontrivial in each degree.
		\item $\widetilde{e}_1^{\ot d}\circ \Delta=\Delta_{\Ain^*}\circ \widetilde{e}_d$ as elements of
		\[\mathrm{Hom}_{{\cal DP}^d_{{\cal F}}}(\Gamma^d(I\otimes I^*),(I^d(I\otimes I^*))_{\Ain^*}),\]
		where $\Delta: \Gamma^d(I\ot I^*)\ra I^{d}(I\ot I^*)$ is the natural embedding and
		\[\widetilde{e}_1^{\ot d}\in \mathrm{Hom}_{{\cal DP}^d_{{\cal F}}}(I^d(I\otimes I^*),(I^d(I\otimes I^*))_{\Ain^*})\]
		is the $d$th external power of  $\widetilde{e}_1$.
	\end{enumerate}
\end{lem}
{\bf Proof of Lemma 3.6}  We obtain our classes from the classes $c[d]^{(i)}$ by applying a multitwist analog of [C3, Lemma 3.4] and then pulling the obtained elements  to 
the category ${\cal F}$. \qed
{\bf Proof of Theorem 3.4}
We will construct $\phi$ in two steps using the intermediate DG category
$\Gamma^d{\cal V}_{X,\Ain}^{op}$. This is yet another category with objects as in ${\cal V}_{\Ain}$ and
\[
\mbox{Hom}^n_{\Gamma^d{\cal V}_{X,\Ain}^{op}}(V,V'):=
\mbox{Hom}_{{\cal F}}(X^{\bullet}(V'\ot \Ain,-),X^{\bullet}(V\ot \Ain,-)[n]).
\]
The composition of morphisms is, like in $\Gamma^d{\cal V}_{X}^{op}$,
 given by the composition in Hom complexes.\newline
% somewhat dual to that in $\Gamma^d{\cal V}_{\Ain}$. For $\alpha\in
%\mbox{Hom}_{\Gamma^d{\cal V}_{X,\Ain}^{op}}(V,V')$ and $\beta\in
%\mbox{Hom}_{\Gamma^d{\cal V}_{X,\Ain}^{op}}(V',V'')$,
%$\beta\circ\alpha$ is given as the composite:
%\[X^{\bullet}(V''\ot\Ain,-)
%\stackrel{m^*}{\ra}
%X^{\bullet}(V''\ot\Ain\ot\Ain,-)\stackrel{\beta_{\Ain}}{\ra}
%X^{\bullet}(V'\ot\Ain,-)\stackrel{\alpha}{\ra}
%X^{\bullet}(V,-),
%\]
%where $m^*$ is induced by the multiplication in $\Ain$ and $\beta_{\Ain}$
%is $\beta$ precomposed with tensoring by $\Ain$ on the first variable.\\
We start with constructing the functor $\rho: 
\Gamma^d{\cal V}_{\Ain}\ra\Gamma^d{\cal V}_{X,\Ain}^{op}$ being the identity on the objects. Thus we need the family of maps:
\[\rho_{V,V'}:\Gamma^d(\mbox{Hom}_{\Ain}(V\ot\Ain,V'\ot\Ain))
\ra\mbox{Hom}_{{\cal F}}(X^{\bullet}(V'\ot\Ain,-),X^{\bullet}(V\ot\Ain,-))\]
compatible with the compositions. For this we observe that since $X^{\bullet}$ is a contravariant strict polynomial with respect to the first variable, 
any $\ka$-linear map (hence in particular any  $\Ain$-linear map) from $V\ot\Ain$ to $V'\ot\Ain$ induces the transformation 
from $X^{\bullet}(V'\ot\Ain,-)$ to $X^{\bullet}(V\ot\Ain,-)$. Since this correspondence is $d$-linear, it amounts to the required linear map
from
$\Gamma^d(\mbox{Hom}(V\ot\Ain,V'\ot\Ain))$ to
$\mbox{Hom}_{{\cal F}}(X^{\bullet}(V'\ot\Ain,-),X^{\bullet}(V\ot\Ain,-))$.
Now the compatibility of $\rho$ with the compositions in the source and target 
categories just boils down to the fact that the action of a functor on the morphisms commutes with their composition. Let us also observe that 
$\rho_{V,V'}$ can be factorized as the composite:
\[
\Gamma^d(\mbox{Hom}_{\Ain}(V\ot\Ain,V'\ot\Ain))\simeq
\Gamma^d(\mbox{Hom}_{\Ain}(V'^*\ot\Ain^*,V^*\ot\Ain^*))\subset
\]
\[
\Gamma^d(\mbox{Hom}(V'^*\ot \Ain^*,V^*\ot\Ain))
\simeq
\mbox{Hom}_{{\cal P}_d}(\Gamma^d(V'^*\ot\Ain^*\ot -),\Gamma^d( V^*\ot\Ain^* -))\simeq
\]
\[
 \mbox{Hom}_{{\cal F}}(\Gamma^d(V'^*\ot \Ain^*\ot -),\Gamma^d( V^*\ot\Ain^*\ot -))\ra\]
\[
\mbox{Hom}_{{\cal F}}(X^{\bullet}(V'\ot \Ain,-),X^{\bullet}(V\ot\Ain,-)),
\]
where the second isomorphism is the Yoneda lemma, the third one follows from the fact that $|\ka|\geq d$. The last arrow is the lifting of morphisms to 
resolutions, hence a priori it exists only up to chain homotopy. Thus our construction of $\rho$ may be thought of as finding such lifts canonically in our situation.\\
Now we turn to constructing the functor 
\[e:\Gamma^d{\cal V}_{X,\Ain}^{op}\ra
\Gamma^d{\cal V}_{X}^{op}  
\]
again being the identity on the objects. 
For this we choose for the element 
 $\widetilde{e}_d\in \mbox{Hom}_{{\cal DP}_{{\cal F}}^d}(\Gamma^d(I^*\otimes I),
 \Gamma^d(I^*\ot \Ain^*\ot I))$  a representing cocycle   
 $\widetilde{e}_d'\in\mbox{Hom}_{{\cal P}_{{\cal F}}^d}(X^{\bullet}(-,-),
 X^{\bullet}(-\ot \Ain,-))$.  
Then precomposing with $\widetilde{e}_d'$
evaluated on the first variable gives
for any $V,V'$ the arrow
\[(\widetilde{e}_d)^*:\mbox{Hom}_{{\cal F}}(X^{\bullet}(V'\ot \Ain,-),X^{\bullet}(V\ot\Ain,-))\ra\mbox{Hom}_{{\cal F}}(X^{\bullet}(V',-),
X^{\bullet}(V\ot\Ain,-)).\]
Next, let 
\[X(i):X^{\bullet}(V\ot\Ain,-)\ra X^{\bullet}(V,-)\]
stand for the transformation induced by the embedding $\ka\subset\Ain$.
Then we define:
\[ e_{V,V'}:=X(i)\circ (\widetilde{e}_d)^*.\]
Again, although the choice of representative $\widetilde{e}_d'$ is not unique,  thanks to its functoriality in the first variable, the precomposition with it (followed by postcomposition with $X(i)$) is a functor from 
$\Gamma^d{\cal V}_{X,\Ain}^{op}$ to
$\Gamma^d{\cal V}_{X}^{op}$.\\
Finally we put $\phi$ to be $e\circ \rho$.\\
It remains to show that $\phi_{V,V'}$ is a quasi-isomorphism for any $V,V'$. The argument follows closely 
the last part of the proof of [C4, Th. 4.2] (which in turn was  a reinterpretation of that of [C3, 
Th. 3.2]), hence we only sketch it. It uses in a crucial way the assumption 
$|\ka|\geq d$, since it relies on Proposition 3.3. We start with the case $d=1$. Here the fact that
$\phi_{V,V'}$ is a quasi-isomorphism follows from the first property of $\widetilde{e}_d$
and Proposition 3.3 (for $d=1$). For a greater $d$ we observe that it follows from the Kunneth formula and the second property of the classes $\widetilde{e}_d$ that $H^*(\phi_{V,V'})$
is onto. Then, since by Proposition 3.3 we know that the domain and codomain of $H^*(\phi_{V,V'})$ have equal graded dimension, the assertion follows. \qed
This formality theorem, which is yet another incarnation of the phenomena observed in
[C3, C4] allows us to perform  ``an $\infty$--affine extension'' of the $\{{\bf f},
{\bf Ra}\}$ adjunction along the lines of [C4, Sect. 3, 4]. The reader is again referred to [C4]
for more extensive explanations of the construction.\newline 
Let ${\cal DP}_X$ be the derived category of the DG category $\mbox{Dif}(\Gamma^d{\cal V}_X^{op})$. Then by
 by Theorem 3.4  we get
an equivalence of triangulated categories
\[{\bf R}\phi^*: {\cal DP}_X\simeq {\cal DP}_d^{af_{\infty}}. \]
Next, it is well known (see e.g. [Ku2]) that ${\cal F}$ may be thought of as the category of
linear functors on  the category  ${\ka{\cal V}}$ (we recall that
$\mbox{Hom}_{\ka{\cal V}}(V,W)=$\\$\ka[\mbox{Hom}_{\ka}(V,W)]$). Thus, in terms of
formalism of [Ke],
$X$ is a $\ka{\cal V}$--$\Gamma^{d}{\cal V}_X^{op}$ bimodule. Hence  we can consider
``the standard functors'' [Ke; C4, Sect. 3]: 
\[H_X: \mbox{Dif}(\ka{\cal V}^{op})\ra \mbox{Dif}(\Gamma^d{\cal V}_X^{op}),\ \ \ \ T_X:\mbox{Dif}(\Gamma^{d}{\cal V}_X^{op})\ra \mbox{Dif}(\ka{\cal V}^{op})\]
and their derived functors
\[{\bf R}H_X: {\cal DF}\ra {\cal DP}_X\ \ \ \ {\bf L}T_X:{\cal DP}_X\ra{\cal DF}.\]
We recall, that since we do not have any boundedness conditions, ${\cal DF}$ stands for
the unbounded derived category.\newline
Now we define ``the $\infty$--affine forgetful functor'':
\[ {\bf f}^{af_{\infty}}:{\cal DP}_{d}^{af_{\infty}}\ra {\cal DF}\]
as ${\bf f}^{af_{\infty}}:={\bf L}T_X\circ ({\bf R}\phi^*)^{-1}$,\newline
and ``the $\infty$--affine  right algebraification'':
\[{\bf a}^{af_{\infty}}:{\cal DF}\ra {\cal DP}_d^{af_{\infty}}\]
as ${\bf a}^{af_{\infty}}:={\bf R}\phi^*\circ {\bf R}H_X$.\newline
The next theorem is the main result of the first part of the paper. It
 is analogous to [C4, Th. 5.1], though slightly weaker. The reason is that, in contrast
 to [C4], $X^{\bullet}$ is not bounded, hence is not a finite object in 
 ${\cal DF}$.
 This forces us to consider the category ${\cal DP}_d^{faf_{\infty}}$ which is the smallest full triangulated subcategory of
${\cal DP}_d^{af_{\infty}}$ containing representable functors and 
closed under isomorphisms and  direct summands. In fact, by [Ke, Th. 5.3] it coincides with the full subcategory of ${\cal DP}_d^{af_{\infty}}$ consisting of finite objects.
Then we have
\begin{theo}
	Functors ${\bf f}^{af_{\infty}}$, ${\bf a}^{af_{\infty}}$ have the following properties:
	\begin{enumerate}
		\item ${\bf f}^{af_{\infty}}\circ z^*\simeq {\bf f}$, $t^*\circ {\bf a}^{af_{\infty}}\simeq {\bf Ra}$.
		\item ${\bf a}^{af_{\infty}}$ is right adjoint to ${\bf f}^{af_{\infty}}$.
		\item ${\bf a}^{af_{\infty}}\circ{\bf f}^{af_{\infty}}\simeq Id_{{\cal DP}_d^{faf_{\infty}}}$.
		\item ${\bf f}^{af_{\infty}}$ restricted to ${\cal DP}_d^{faf_{\infty}}$ is fully faithful.
		%\item The triangulated quotient category ${\cal DP}_{pd}/{\cal DP}_d^{af}$  is %equivalent to the Verdier localization of ${\cal DP}_{pd} $ with respect to the %essential image of ${\bf C}^{af}$ (see e.g. [Kr]).
	\end{enumerate}
\end{theo}
{\bf Proof:} 
In order to get the first isomorphism in the first part we evaluate 
${\bf f}^{af_{\infty}}\circ z^*$ on the projective generator $\Gamma^{d,U}$ of ${\cal P}_d$. We obtain
\[{\bf f}^{af_{\infty}}\circ z^*(\Gamma^{d,U})={\bf f}^{af_{\infty}}(h^{U\ot A_{\infty}})=X(U,-)\simeq \Gamma^{d,U}={\bf f}(\Gamma^{d,U}),\]
hence we get a natural in $V$ isomorphism ${\bf f}^{af_{\infty}}\circ z^*(\Gamma^{d,U})\simeq
{\bf f}(\Gamma^{d,U})$. This, since any object in ${\cal DP}_d$ can be represented by a complex of coproducts of $\Gamma^{d,U}$ and the both functor commute with infinite coproducts, gives the first isomorphism. To get the second isomorphism we observe that
\[t^*\circ {\bf a}^{af_{\infty}}(F)(V)=\mbox{Hom}_{{\cal F}}(X(V,-),F)={\bf Ra}(F)(V)\]
for any $F\in{\cal DF}$.\newline
The second part of the theorem follows from the $\{{\bf L}T_X,{\bf R}H_X\}$ adjunction and the fact that ${\bf R}\phi^*$ is an equivalence.\newline
To get the third part  of the theorem we first observe that \[{\bf a}^{af_{\infty}}\circ{\bf f}^{af_{\infty}}(h^{U\ot\Ain})\simeq h^{U\ot \Ain}\]
by the very definition of $X$. From this we conclude that the unit of the adjunction
is an isomorphism on the whole category  ${\cal DP}_d^{faf_{\infty}}$.\newline
Part 4   follows formally  from  parts 2, 3 (by e.g. [Kr, Prop.~2.3.1]).\qed
\section{Spectra of strict polynomial functors}
In this section we modify the category (of complexes over) ${\cal P}$ by formally inverting the Frobenius
twist operation.  We achieve this goal by a general construction, known from stable homotopy theory, i.e. we consider spectra of complexes of strict polynomial functors.
We follow a general approach of Hovey [Ho1]  who starts from a Quillen model category ${\cal C}$ with a left
Quillen endofunctor $T$ and  equips its  category of spectra with an appropriate 
Quillen model structure. In order to conform to this context we have to slightly
adjust our setup. Namely, although we are mainly interested in ${\cal P}_d$ for a fixed $d>0$, we should also allow  strict polynomial functors of other degrees
to make the Frobenius twisting functor ${\bf C}$ into an endofunctor (of course it would suffice
to consider degrees $dp^i$). However, the category ${\cal P}:=\bigoplus_{d>0}{\cal P}_d$ is not suitable for our purposes, since the forgetful functor ${\bf f}:
{\cal P}\ra {\cal F}$ does not possess the right adjoint. Hence we shall consider
the product category
\[\widehat{{\cal P}}:=\prod_{d>0}{\cal P}_d\]
and we recall again that we do not assume that the objects in ${\cal P}_d$
are functors taking finite dimensional values. Then we consider
the category ${\cal K}\widehat{{\cal P}}$ of unbounded complexes over
$\widehat{{\cal P}}$. ${\cal K}\widehat{{\cal P}}$ can be equipped with the projective Quillen model structure. For the readers convenience we recall this structure.
It is described in detail e.g. in
[Ho2] for the category of modules over a ring but this description readily generalizes to any Grothendieck category ${\cal A}$. In this structure all complexes are fibrant, while the cofibrant ones are those satisfying the ``property P'' [Ke]. In our situation
the property P boils down to saying that a complex admits a filtration 
$\{M_j\}_{j\geq 0}$ such that for any $j\geq 0$:
\begin{itemize}
\item
embedding $M_j\subset M_{j+1}$ splits in ${\cal A}$.
\item 
$M_{j+1}/M_j$ consists of projectives and has trivial differential.
\end{itemize}
The fibrations are the epimorphisms, the cofibrations are the  monomorphims with cofibrant cokernels which split in the underlying abelian category. The weak equivalences are the quasi-isomorphisms. Then we encounter another technical problem: since ${\bf C}$ does not preserve projective objects, it is not a left Quillen functor with respect to the projective Quillen structure
on  ${\cal K}\widehat{{\cal P}}$. The simplest way of overcoming this obstacle is by using the following technical fact.
\begin{prop}
There exists a functor ${\bf C}': 	{\cal K}\pec\ra {\cal K}\pec$ such that
\begin{enumerate}
	\item ${\bf C}'$ is left Quillen functor with respect to the projective model structure on ${\cal K}\pec$.
\item There is a natural isomorphism of total left derived functors
 ${\bf LC}'\simeq {\bf LC}$.
\end{enumerate}
\end{prop}
{\bf Proof:} We can describe ${\bf C}'$ by using an explicit construction from [C4].
Namely we established  in [C4, Theorem 5.1 (1)]  an isomorphism of functors 
${\bf C}\simeq{\bf C}^{af}\circ z^*$ as functors from ${\cal DP}_{dp}$ to ${\cal DP}_d$.
We emphasize that $z^*$ stands here for the functor appearing in [C4]  (we also remind that  we stick to the convention that for exact functors we denote their derived functors by the same letter). 
Then  we recall from  [C4, Sections 4,5]
that ${\bf C}^{af}$ is the derived functor of the composite $T_X\circ (\phi^*)^{-1}$
of two left Quillen functors.  Hence, since readily $z^*$ is left Quillen functor,
 we can define ${\bf C}'$ as the composite  $T_X\circ \phi^{-1}\circ z^*$. 
 In order to obtain an equivalence ${\bf LC}'\simeq {\bf LC}$ we invoke again 
 [C4, Theorem 5.1 (1)]. Namely, 
 in the proof of [C4, Theorem 5.1] we constructed a collection of quasi-isomorphisms
\[{\bf C}'(\Gamma^{d,U})\ra {\bf C}(\Gamma^{d,U})\]
natural in $U$. This allows us 
to construct a transformation ${\bf C}'\ra {\bf C}$ of functors defined on the full subcategory of ${\cal KP}_d$ consisting of complexes projective in each degree.
Since the cofibrant objects in ${\cal KP}_d$ are projective in each degree, we obtain
a natural isomorphism ${\bf LC}'\simeq {\bf LC}$. Finally  we extend our construction degreewise to the whole product category ${\cal K}\pec$. \qed
Let ${\bf K}': {\cal K}\pec\ra {\cal K}\pec$ be the right adjoint functor to ${\bf C}'$.
From now on we consider the pair of adjoint functors $\{{\bf C}',{\bf K}'\}$ instead of the original adjoint pair $\{{\bf C},{\bf K}\}$. 
In fact, we will rarely refer to the specific construction of ${\bf C}'$ given above.
In most cases, the properties listed in Proposition 4.1 will be sufficient for our
purposes.
\newline
Now we apply the machinery of [Ho1]
to the category ${\cal K}\widehat{{\cal P}}$ equipped with the projective model structure with the endofunctor ${\bf C}'$. Namely,  we form the category of spectra over ${\cal K}\pec$.
\begin{defi}
	We call a collection of complexes $F_i\in {\cal K}\pec$ and cochain maps $\tau_i: {\bf C}'(F_{i})\ra
	F_{i+1}$ for all $i\geq 0$ a spectrum (of complexes of strict polynomial functors). For spectra 
	$F_{\bullet}, G_{\bullet}$ we call a collection of cochain maps $\phi_i: F_i\ra G_i$
	a map of spectra if $\tau_i\circ {\bf C}'(\phi_i)=\phi_{i+1}\circ \tau_i$ for all $i\geq 0$.
\end{defi}
Readily the spectra and maps of spectra inherit from ${\cal K}\pec$ the structure of DG category (ie. grading comes from the shift in ${\cal K}\widehat{{\cal P}}$ and not that in spectra). We call 
this category the category of
spectra (of complexes of strict polynomial functors) and 
denote it by ${\cal S}\pec$.

For any $F\in{\cal K}\pec$ we have a spectrum ${\bf C}^{\infty}(F)$ defined by the formula
\[{\bf C}^{\infty}(F)_i:={\bf C}'^i(F),\]
where ${\bf C}'^i$ is the $i$th iteration of ${\bf C}'$.
The assignment $F\mapsto {\bf C}^{\infty}(F)$ produces the functor 
\[{\bf C}^{\infty}:{\cal K}{\pec}\ra {\cal S}\pec\]
which has the evaluation functor $ev(F_{\bullet}):=F_0$ as right adjoint.
\begin{defi}
We call a spectrum $F_{\bullet}$ a ${\bf C}$--spectrum  if all the maps 
$\tau_i: {\bf C}'(F_i)\ra F_{i+1}$ are quasi-isomorphisms. 
Similarly, we call a spectrum $F_{\bullet}$ a ${\bf K}$--spectrum if all the maps 
$\omega_i: F_i\ra {\bf K}'(F_{i+1})$ adjoint to $\tau_i$ are quasi--isomorphisms. 
\end{defi}
Of course, ${\bf C}^{\infty}(F)$ is always a  ${\bf C}$--spectrum. Less trivially,
for a Young diagram $\la$, let $p\la$ denote the Young diagram whose rows are those
of $\la$ multiplied by $p$. Then  by [C3, Prop. 2.1], ${\bf K}'(S^{p\la}) \simeq S^{\la}$, hence the collection $\{ S^{p^i\la}\}_{i\geq 0}$
forms a  ${\bf K}$--spectrum denoted $S^{p^{\bu}\la}$. More generally, since by [C3, Prop. 4.2],
${\bf K}(S_{F_k(\la)}[h_k])=S_{\la}$, any Schur functor $S_{\la}$ gives rise to a  
 ${\bf K}$--spectrum $\{S_{F_k^i(\la)}[h_k^i]\}$ ($F^i_k$ stands for certain combinatorial operation which enlarges Young diagram (see [C2])).
 In fact, any spectrum can be turned into a ${\bf K}$--spectrum by means of
 the ``delooping functor'':
 \[\Theta^{\infty}: {\cal S}\pec\ra {\cal S}\pec\] 
  given by the formula:
 \[\Theta^{\infty}(F_{\bullet})_i:=\mbox{colim}_j {\bf K}'^j(F_{i+j}).\] 
 Then, as it is explained in [Ho1, Section 1], for any model category with left Quillen endofunctor one can endow its category of spectra with the obvious model structure, which is (somewhat unfortunately) also called ``projective''. In order to avoid confusion
with our projective model structure on ${\cal K}\pec$ we shall call this model structure ``levelwise''. Our terminology is justified by the fact that a morphism of spectra $\phi_{\bullet}: F_{\bullet}\ra G_{\bullet}$ is a cofibration (resp. a weak equivalence) if and only if all $\phi_i$ are cofibrations (resp. weak equivalences). The class of fibrations is determined by the lifting 
property.
 Then the final model structure on the category of spectra, called in [Ho1]  ``the stable model structure'' is obtained from the levelwise model structure  by the Bousfield localization process with respect to certain class of morphisms (see
[Ho1, Sect. 2]). We summarize below the basic properties of the stable model 
 Quillen structure on ${\cal S}\pec$.
\begin{prop}
There exists a finitely generated model structure on ${\cal S}{\pec}$ with the following properties:
\begin{enumerate}
\item The pair of functors  $\{{\bf C}^{\infty},ev\}$ is a Quillen pair. 
\item The cofibrant objects are the spectra consisting of  complexes 
satisfying the property P
with structure maps $\tau_i$ monomorphic with cokernels satisfying the property P.
\item The fibrant objects are the ${\bf K}$-spectra.
\item A map of spectra which is a levelwise quasi-isomorphism is 
a weak equivalence.
\item The degreewise prolongation of ${\bf C}'$ on ${\cal S}\pec$ is a Quillen equivalence with the quasi-inverse  being the shift functor.
%\item If for a map of spectra $\phi_{\bu}: F_{\bu}\ra G_{\bu}$ there exists $i_0$ %such
%that for all $i\geq i_0$ $\phi_i$ is a quasiisomorphism, then $\phi_{\bu}$ is a %weak
%equivalence.
%\item The pair of functors  $\{{\bf C}^{\infty},ev\}$ is a Quillen pair. 
\item A natural map $X\ra \Theta^{\infty}(X)$ is a weak equivalence for any
spectrum $X$, hence $\Theta^{\infty}$ can be chosen as a fibrant replacement functor.
\item
The category  $Ho({\cal S}{\pec})$  has a structure of triangulated category
such that the total derived functor ${\bf L}{\bf C}^{\infty}={\bf C}^{\infty}: Ho({\cal S}\pec)\ra
Ho({\cal S}{\pec})$ is an exact functor.

%\item For any $F: X\ra Y$ with $X,Y\in {\cal SP}_d$ and $d\in\np$, its functorial %factorizations
%into fibrations/cofibrations belong to ${\cal SP}_d$. Hence ${\cal SP}_d$ has a %natural, inherited from ${\cal SP}$ model structure. Moreover, we have %decomposition of homotopy categories:
%\[Ho({\cal SP}\simeq \prod_{d\in\np}Ho({\cal SP}_d).\]
%\item The prolongation of ${\bf C}'$ restricts to a Quillen equivalence 
%${\cal SP}_d\ra {\cal SP}_{pd}$ for any $d\in\np$.

\end{enumerate}	
\end{prop}
{\bf Remark:  } These properties of the model structure on ${\cal S}{\pec}$ will be sufficient
for our purposes. The first property allows one to compare ${\cal K}{\pec}$ with ${\cal S}\pec$. 
Properties 2 and 3 provide explicit descriptions of fibrant and cofibrant objects which is essential in making  calculations and also shows importance of ${\bf K}$-spectra. In fact one could also extract from [Ho1, Sections 3-4] certain
descriptions of fibrations and cofibrations in ${\cal S}{\pec}$ but we will not need them.
Properties 4 and 5 deal with weak equivalences. Property 5 is central for the whole 
idea of spectra, since it shows that in ${\cal S}{\pec}$  the endofunctor ${\bf C}'$
becomes invertible up to homotopy. Property 6 will be crucial in the proof of
Theorem 4.6.\newline
{\bf Proof: } We apply the machinery of [Ho1] to the projective model structure on
${\cal K}\pec$ with the functor ${\bf C}'$ as the left Quillen endofunctor. 
The first property is  the second assertion in [Ho1, Proposition 1.15].
The second property 
follows from [Ho1, Proposition 1.14] (we recall that the Bousfield localization preserves cofibrant objects by [Ho1, Theorem 2.2]). The third part is [Ho1, Theorem 3.4].
The fourth property is obvious.
The fifth property  is [Ho1, Theorem 3.9].\newline
We would like to deduce the sixth property from [Ho1, Corollary 4.11]. For this we need to show that sequential colimits in ${\cal K}\pec$ commute with finite products and that ${\bf K}'$ commutes with sequential colimits. The first fact is well known, since  ${\cal P}_d$ is equivalent to a module category.  In order 
to show the second fact it suffices to show that 
 ${\bf K}'$, when restricted to ${\cal KP}_{dp}$ commutes with sequential colimits. It follows from [C4, Section 5] that
 ${\bf K}'$ is  explicitly given as
 \[{\bf K}'(F)(V)=\mbox{Hom}_{{\cal P}_{dp}}(\widetilde{X}(V,-),F),\]
 where $\widetilde{X}$ is certain strict polynomial functor in two variables. Then, since 
 $\widetilde{X}$
 is finite dimensional, ${\bf K}'$ commutes with infinite sums; since $\widetilde{X}(V,-)$ consists of projectives in ${\cal P}_{dp}$, ${\bf K}'$ preserves cokernels. Therefore
 ${\bf K}'$ preserves sequential colimits.\newline  In order to equip $Ho({\cal S}{\pec})$ with
 a structure of triangulated category we recall that by [Ho2, Sect. 7] 
 the homotopy category $Ho({\cal C})$  of a Quillen model category 
 ${\cal C}$
 is canonically a triangulated
 category  whenever  the ``model theoretic suspension functor'' on ${\cal C}$
 [Ho2, Sect. 6] becomes an equivalence on $Ho({\cal C})$. The standard construction  of the triangulated structure on the derived category of abelian category 
 fits into this formalism. Hence we interpret the standard triangulated structure on ${\cal D}\pec$ with the shift of complexes as the suspension  in model theoretic terms.  Then we prolong degreewise the shift functor on ${\cal S}{\pec}$ and, since
 it remains  invertible (already on ${\cal K}{\pec}$), this allows one to make
 $Ho({\cal S}{\pec})$ with degreewise model structure into  a triangulated category.   
 At last, since the
 Bousfield localization preserves the cofibrations, it commutes with the suspension functor. This shows that the suspension functor remains a weak equivalence in 
 the stable model structure on ${\cal S}{\pec}$.  
 %In order to obtain the seventh property 
 %we consider the idempotent functor $\pi:{\cal SP}\ra{\cal SP}$ given by the $i$th
 %level by the projection ${\cal K}\pec\ra {\cal P}_{p^id}\subset {\cal K}\pec$.
 %Then, obviously $im(\pi)={\cal SP}_d$ and it suffices to show that $\pi$
 %preserves the classes of fibrations, cofibratons and weak equivalences. 
 %We first  observe that the analogous facts are obvious for  the degreewise %model structure on ${\cal SP}$, since ${\cal SP}=\prod_{d\in\np}{\cal SP}_d$.  %The resulting ``stable'' model structure on ${\cal SP}$ is obtained by the %Bousfield  
%localization.  Since the Bousfield localization does not change the class of %cofibrations [Ho,?] we see that $\pi$ preserves the cofibratons. Now we look at %the weak equivalences. This is  the class of ``local maps'' with respect to the %``generating class'' of maps of the form  $S\circ {\bf C}'(\Sigma^{\infty}(X))\ra %\Sigma^{\infty}(X)$
%for $X\in{\cal S}\pec$, where $S(X_{\bullet}):=X_{\bullet-1}$ is the shift functor.
%Since $\pi(\Sigma^{\infty}(X))=\Sigma^{\infty}(\pi(X))$, we see that $\pi$ %preserves the generating class, hence the class of local maps i.e. the class of
%weak equivalences. The fact that $\pi$ preserves the class of fibrations follows
%from the fact that it is determined by the lifting property. The eight property
%follows from properties 5 and 7. 
\qed
The last property justifies calling  $Ho({\cal S}{\pec})$ the derived category of ${\cal S}{\pec}$, hence from now on we shall use notation ${\cal DS}{\pec}:=Ho({\cal S}{\pec})$.\newline
Now we would like to find inside ${\cal DS}{\pec}$ a subcategory corresponding to
${\cal DP}_d$. Of course we have a decomposition of ${\cal K}\pec$ into the product 
of ${\cal KP}_d$. However, in order to describe the situation for spectra  we need more subcategories.
Let ${\bf N}[\frac{1}{p}]:=\{ep^j: e,j\in{\bf Z}, e>0, (e,p)=1\}$. Then 
for $e p^j\in {\bf N}[\frac{1}{p}]$ 
we define ${\cal SP}_{ep^j}$
to be the full subcategory of ${\cal SP}$ consisting of spectra $X_{\bullet}$
such that $X_i=0$ for $i<-j$ and $X_i\in {\cal KP}_{ep^{i+j}}$ for $i\geq -j$.
Then we have a decomposition of DG categories:
\[{\cal S}{\pec}=\prod_{d\in{\bf N}[\frac{1}{p}]}{\cal SP}_d.\] 
Now we claim that, roughly speaking, all our constructions carry over to the subcategories ${\cal SP}_d$. We recall that we could not apply the ideas of [Ho1]
to ${\cal SP}_d$ directly, since ${\bf C}'$ does not preserve these subcategories.
\begin{prop}

	For any $f: X\ra Y$ with $X,Y\in {\cal SP}_d$ and $d\in\np$, its functorial factorizations
	into fibrations/cofibrations belong to ${\cal SP}_d$. Hence ${\cal SP}_d$ has a natural
	model structure  inherited from ${\cal SP}$. 
	The homotopy category ${\cal DSP}_d:=Ho({\cal SP}_d)$ has a structure of triangulated category and there is an equivalence of triangulated categories:
	\[{\cal DS}{\pec}\simeq \prod_{d\in\np}{\cal DSP}_d.\]
	Moreover, the prolongation of ${\bf C}'$ restricts to a Quillen equivalence 
	producing an equivalence of triangulated categories
	\[{\cal DSP}_d\simeq {\cal DSP}_{pd}\] for any $d\in\np$.

\end{prop}
{\bf Proof: }
For $d=d'p^i\in\np$ we consider the  functor $\pi_d:{\cal S}{\pec}\ra{\cal S}{\pec}$ given on the $j$th
level by the projection ${\cal K}\pec\ra {\cal KP}_{p^{j-i}d}\subset {\cal K}\pec$
if $j\geq i$ and trivial elsewhere.
Then, obviously $im(\pi_d)={\cal SP}_d$ and $\pi_d$ restricted to ${\cal SP}_d$ is
the identity functor. Now, since for any $f\in \mbox{Mor}({\cal S}{\pec})$, $\pi_d(f)$
is a retract of $f$, $\pi_d$ preserves the classes of cofibrations, fibrations and weak equivalences. This shows that ${\cal SP}_d$ has a natural model structure:
we obtain a functorial factorization into cofibrations/fibrations in ${\cal SP}_d$
by applying $\pi_d$ to the factorization in ${\cal S}{\pec}$. This also shows that the
decomposition
${\cal S}{\pec}=\prod_{d\in{\bf N}[\frac{1}{p}]}{\cal SP}_d$ is an isomorphism of model categories where we consider on the right hand side the product of model structures we have just described. This gives the required decomposition of homotopy categories, which also preserves triangulated structure because
$\pi_d$ commutes with the cone functor.\newline At last, we observe that by
Proposition 4.4.(5) the prolongation of ${\bf C}'$ is a self--equivalence of
${\cal DS}{\pec}$ which takes ${\cal DSP}_d$ into ${\cal DSP}_{pd}$ which proves the last part of the proposition. \qed

The properties of the derived category of spectra we have established so far allow
us to obtain an analog of the known  description of stable homotopy maps between suspension spectra.  This theorem  is one of the main objectives of this part of the article.

\begin{theo}
 	Let $F\in {\cal KP}_d$ be finite dimensional and $G_{\bullet}\in{\cal SP}_d$. 
 	Then there is 
 	a natural in $F,G_{\bullet}$
 	isomorphism
 	\[{\rm Hom}_{{\cal DSP}_d}({\bf C}^{\infty}(F),G_{\bullet})\simeq
 	{\rm colim}_i {\rm \ \! Hom}_{{\cal DP}_d}(F,{\bf K}'^i(G_{i})).\]
 	In particular, for  $F,G\in {\cal P}_d$ we obtain:
 	\[{\rm Hom}_{{\cal DSP}_d}({\bf C}^{\infty}(F),{\bf C}^{\infty}(G)[s])\simeq
 	{\rm colim}_i {\ \!\rm Ext}^s_{\pu_{dp^i}}(F^{(i)},G^{(i)}).\]	
 \end{theo}
 {\bf Proof:} We shall deduce our theorem from [Ho1, Corollary 4.13], hence we should verify the assumptions of [Ho1, Corollary 4.13]. Let 
 $A$ be a cofibrant replacement of $F\in {\cal KP}_d$. Then $A$ and its cylinder
 are finite objects in ${\cal K}\widehat{{\cal P}}$ as being finite dimensional complexes of projectives. Additionally, 
 we recall that any spectrum $G_{\bullet}$ is  fibrant in the levelwise model structure. Therefore we can 
 apply [Ho1, Corollary 4.13] to our $A$ and $Y:=G_{\bullet}$ and we obtain   
 the first part of our theorem. The second part is just a special case which we
 distinguished since it is related to [FFSS]. Indeed, we obtain:
 \[{\rm Hom}_{{\cal DSP}_d}({\bf C}^{\infty}(F),{\bf C}^{\infty}(G)[s])\simeq
 \mbox{colim}_i {\rm \ \! Hom}_{{\cal DP}_d}(F,{\bf K}'^i({\bf C}'^i(G)[s])))\simeq\]
 \[\mbox{colim}_i {\rm \ \! Hom}_{{\cal DP}_{p^id}}({\bf C}'^i(F),{\bf C}'^i(G)[s])
 \simeq {\rm colim}_i {\ \!\rm Ext}^s_{\pu_{dp^i}}(F^{(i)},G^{(i)}).\]	
 \qed
In fact, the crucial ingredient of the proof of [Ho1, Corollary 4.13] is the fact that
the delooping functor $\Theta^{\infty}$ is a weak equivalence. It is worth mentioning (we will use this fact in Section 6) that in our situation,
thanks to the Collapsing Conjecture [C3, Theorem 3.2] we have a very explicit description of the functor $\Theta^{\infty}{\bf C}^{\infty}$. Namely,
 for $F\in {\cal K}\pec$, let $F^{(\infty)}_{A_{\infty}}$ denote the spectrum with
 \[(F^{(\infty)}_{A_{\infty}})_i:={\bf C}'^i(F)_{A_{\infty}},\]
 where ${\bf C}'^i(F)_{A_{\infty}}$ is meant as precomposing of ${\bf C}'^i(F)$ with $-\ot A_{\infty}$ which practically means that we also  twist $A_{\infty}$. This  is of great importance because  twisting a
 graded space multiplies degrees of elements by $p^i$. 
 When defining the map $\tau_i: {\bf C}'^{i}(F)_{A_{\infty}}\ra {\bf C}'^{i+1}(F)_{A_{\infty}}$ 
 we should also be careful since 
 \[{\bf C}((F^{(i)})_{A_{\infty}})(V)=F(V^{(i+1)}\ot A_{\infty}^{(i)})\neq 
 F(V^{(i+1)}\ot A_{\infty}^{(i+1)})=(F^{(i+1)})_{A_{\infty}}(V)\]
 Using these descriptions, we take $\tau_i$ as the map induced by the projection
 $A_{\infty}^{(i)}\ra A_{\infty}^{(i+1)}$. Thus we see that 
 $\tau_i$ is not an isomorphism.  This is the reason why we do not use  the notation
 ${\bf C}^{\infty}(F_{A_{\infty}})$ here. On the other hand, $F^{(\infty)}_{A_{\infty}}$
 is a  ${\bf K}$--spectrum, since ${\bf K}((F^{(i+1)})_{A_{\infty}})=(F^{(i)})_{A_{\infty}}$
 and, as it is easy to see, $\omega_i$ corresponds just to the identity map.
 Now we have 
 \begin{prop}
  Let $\{F_{\bu}\}$ be a ${\bf C}$--spectrum.
 	Then there is a  natural in $F_{\bu}$ weak equivalence
 \[(F_0)^{(\infty)}_{A_{\infty}}\simeq 
 \Theta^{\infty}(F_{\bu}).\]
 \end{prop}
 {\bf Proof: }
 By the Collapsing Conjecture [C3, Theorem 3.2] we have
 \[\Theta^{\infty}(F_{\bu})_i\simeq{\rm colim}_j {\bf K}^j(F_0^{(i+j)})\simeq
 {\rm colim}_j (F_0^{(i)})_{A_j}=(F_0^{(i)})_{A_{\infty}}=((F_0)^{(\infty)}_{A_{\infty}})_i.\]
 \qed
 \section{Spectra of ordinary functors and factorization}
The aim of the present section is to factorize the adjunction $\{{\bf f}, {\bf Ra}\}$
through ${\cal DSP}_d$. In order to compare the categories ${\cal DSP}_d$ and ${\cal DF}$  we take the following strategy. We introduce an intermediate category ${\cal SF}$ of spectra of ordinary functors. 
Since the Frobenius twist on ${\cal F}$ is invertible,  ${\cal DSF}$
and ${\cal DF}$ are equivalent.
Then we compare ${\cal DSP}_d$ and ${\cal DSF}$ by using functoriality of the construction of spectra.\newline
Let ${\cal KF}$ be the category of complexes of objects of ${\cal F}$ (we admit unbounded complexes). The Frobenius twist 
\[{\bf C}: {\cal KF}\ra {\cal KF}\]
is a self--equivalence, hence a left Quillen endofunctor for the projective model
structure on  ${\cal KF}$. Then we introduce the category ${\cal SF}$ of spectra
over ${\cal KF}$ and we equip it with the stable model structure by the Hovey construction analogous to that applied in Section 4 to the category ${\cal K}\pec$.
Now, since ${\bf C}$ is a Quillen equivalence, we have:
\begin{prop}
The Quillen pair $\{{\bf C}^{\infty},ev\}$ is a Quillen equivalence between ${\cal KF}$
and ${\cal SF}$ (i.e. their derived
functors are mutually inverse equivalences between derived categories).
Explicitly, we can take ${\bf LC}^{\infty}={\bf C}^{\infty}$ and
${\bf R}ev(F_{\bullet})=\mbox{colim}_j F_{nj}$ where $|{\bf k}|=p^n$.
\end{prop}
{\bf Proof: } The fact that we have a Quillen equivalence follows from [Ho1, Theorem 5.1]. The fact that ${\bf LC}^{\infty}={\bf C}^{\infty}$ is a consequence of
the fact that ${\bf C}^{\infty}$ preserves all weak equivalences. Then, in general,
${\bf R}ev=ev\circ R$ where $R$ is a fibrant replacement functor. 
But we observe that in the category ${\cal F}$, the right adjoint functor to ${\bf C}$ is
just ${\bf C}^{-1}={\bf C}^{n-1}$.  Therefore, obviously
$\Theta^{\infty}(F)\simeq F$ for any $F\in {\cal SF}$ and we can take $\Theta^{\infty}$ as $R$. Finally, we have 
\[ev\circ\Theta^{\infty}(F_{\bullet})=\mbox{colim}_j 
{\bf C}^{nj-j}(F_{j})\simeq\mbox{colim}_j F_{nj}.\]
\qed
Now we are going to compare the categories ${\cal DS}\pec$ and ${\cal DF}$.
We introduce a temporary notation.
Let us denote by ${\bf a}_d: {\cal KF}\ra{\cal KP}_d$ the right adjoint functor 
to the forgetful functor ${\bf f}: {\cal KP}_d\ra{\cal KF}$. Then it is easy to see
that the functor 
\[\widehat{{\bf a}}: {\cal KF}\ra{\cal K}\pec\]
given as the product $\widehat{{\bf a}}(F):=({\bf a}_1,\ldots,{\bf a}_d,\ldots)$ is right adjoint
to the functor 
\[\widehat{{\bf f}}: {\cal K}\pec\ra{\cal KF}\]
sending $(F_1,\ldots,F_d,\dots)\in {\cal K}\pec$ to the direct sum $\bigoplus_{d\geq 1}
{\bf f}(F_d)$. Moreover, by Proposition 4.1, we have a natural transformation 
\[\widehat{{\bf f}}\circ{\bf C}'\ra{\bf C}\circ \widehat{{\bf f}}\]
of functors from  ${\cal K}\pec$ to  ${\cal KF}$, 
which is a weak equivalence . Therefore, by [Ho1, Proposition 5.5] we get the Quillen pair
$\{S\widehat{{\bf f}},S\widehat{{\bf a}}\}$ between the spectra categories ${\cal SP}$ and ${\cal SF}$
such that 
\[S\widehat{{\bf f}}\circ {\bf C}^{\infty}\simeq {\bf C}^{\infty}\circ \widehat{{\bf f}}.\] 
Now we are ready for defining the adjunction between ${\cal D}\pec$ and ${\cal DF}$. Let 
\[{\bf f}^{st}:{\cal DS}\pec\ra {\cal DSF}\]
be given as $ev\circ\Theta^{\infty}\circ {\bf L}S\widehat{{\bf f}}$ and
 \[{\bf a}^{st}:{\cal DSF}\ra{\cal DS}\pec\]
given as ${\bf R}S\widehat{{\bf a}}\circ {\bf C}^{\infty}$.
Let ${\cal DP}^{b}_d$ be the 
the full subcategory of ${\cal DP}_d$ consisting of finite dimensional complexes
and let ${\cal DP}^{st}_d$ be the 
smallest full triangulated subcategory of ${\cal DP}_d$ containing
${\bf C}^{\infty}({\cal DP}^{b}_d)$ and closed under isomorphisms and direct summands.
Our terminology here refers to stable homotopy theory where the stable category of Spanier-Whitehead can be characterized in a similar manner as a full subcategory of the category of spectra.
Then we have
 \begin{theo}
 The functors ${\bf f}^{st}$ and ${\bf a}^{st}$ satisfy the following properties:
 \begin{enumerate}
 \item The functor ${\bf f}^{st}$ is left adjoint to
 ${\bf a}^{st}$. 
 \item There are isomorphisms of functors between ${\cal D}\pec$ and ${\cal DF}$
  \[{\bf f}^{st}\circ {\bf C}^{\infty}\simeq \widehat{{\bf f}},\]
  \[ {\bf R}ev\circ {\bf a}^{st}\simeq{\bf R\widehat{a}}.\]
 \item Let $d\leq |{\bf k}|$. Then
 the unit map $Id\longrightarrow {\bf a}^{st}\circ{\bf f}^{st}$ is an isomorphism on the subcategory ${\cal DP}^{st}_d$.
 \item  Let $d\leq |{\bf k}|$. Then ${\bf f}^{st}$ restricted to ${\cal DP}^{st}_d$ is fully faithful.
 \end{enumerate}
\end{theo}
 {\bf Proof:}  The first part follows from the facts that $\{S{\bf f},S{\bf a}\}$ is a Quillen pair and Proposition 5.1. For the second part we recall that 
 \[S\widehat{{\bf f}}\circ {\bf C}^{\infty}\simeq {\bf C}^{\infty}\circ\widehat{ {\bf f}}.\] 
 Therefore we obtain
 \[{\bf f}^{st}\circ {\bf C}^{\infty}\simeq 
 {\bf R}ev\circ S\widehat{{\bf f}}\circ {\bf C}^{\infty}\simeq 
 {\bf R}ev\circ {\bf C}^{\infty}\circ \widehat{{\bf f}}\simeq \widehat{{\bf f}}.\]
 In order to  obtain the second isomorphism we recall, that since $S\widehat{{\bf a}}$ is just
 degreewise prolongation of ${\bf a}$ by [Ho1, Lemma 5.3], it commutes with
 ${\bf C}^{\infty}$. Hence we get
 \[{\bf R}ev\circ {\bf a}^{st}\simeq {\bf R}ev\circ  S\widehat{{\bf a}}\circ{\bf C}^{\infty}
 \simeq {\bf R}ev\circ {\bf C}^{\infty}\circ \widehat{{\bf a}}\simeq{\bf R\widehat{a}}.\]
 Parts 3 and 4  are equivalent. Since ${\cal DP}_d^{st}$ is generated as triangulated category with direct summands by ${\bf C}^{\infty}({\cal DP}^{b}_d)$, it suffices 
 to show that the functor ${\bf f}^{st}$ restricted to the objects of the form ${\bf C}^{\infty}(F)$
 for $F\in {\cal DP}^{st}_d$ is fully faithful.
 To this end
  let us take $F,G\in {\cal P}_d$. 
 Then by Theorem  4.6
 \[{\rm Hom}_{{\cal DSP}_d}({\bf C}^{\infty}(F),{\bf C}^{\infty}(G)[s])\simeq
 {\rm colim}_i {\ \!\rm Ext}^s_{\pu_{dp^i}}(F^{(i)},G^{(i)}).\]	
 On the other hand 
 \[\mbox{Hom}_{{\cal DSF}}({\bf f}^{st}({\bf C}^{\infty}(F)),{\bf f}^{st}({\bf C}^{\infty}(G)[s]))\simeq \mbox{Hom}_{{\cal DSF}}({\bf f}(F),{\bf f}(G)[s])\simeq\mbox{Ext}^s_{{\cal F}}(F,G)
.\]	
Therefore our assertion follows from [FFSS, Theorem 3.10]. \qed
 {\bf Remark 5.3} The functor ${\bf f}^{st}$ has an intriguing extra feature. 
 We recall from Section 4 the fibrant spectra $S^{p^{\bu}\la}:=\{S^{p^i\la}\}$.
 Then we see that
 \[{\bf f}^{st}(S^{p^{\bu}\la})={\rm colim}_i\ S^{p^i\la(-i)}\]
 which is nothing but the product of the Carlsson functors whose injectivity was
 shown by Kuhn [Ku1]. Thus we see that ${\bf f}^{st}$ preserves some important
 fibrant objects in contrast to the fact that the original forgetful functor 
 ${\bf f}: {\cal DP}_d\ra{\cal F}$ does not
 preserve injectives. This suggests possibility of fully reconstructing ${\cal DF}$
 from some categories of algebrogeometric origin. We hope to extend this observation in a future work.
\section{Comparison of spectra and $\infty$--affine functors}
In this section we  construct a functor $\gamma: {\cal DP}_d^{af_{\infty}}\ra {\cal DSP}_d$ which is a  full embedding and is compatible with our previous constructions.\newline
Let $A_i:=\ka[x_1, x_2,\ldots,x_i]/(x_1^p, x_2^p,\ldots x_i^p)$ for $|x_j|=2p^j$
(we allow here also $A_0:=\ka$).
We consider the categories $\Gamma^d{\cal V}_{A_i}$, ${\cal P}_d^{af_i}$, ${\cal DP}_d^{af_i}$, analogous to the notions introduced in Section 2. 
n particular ${\cal P}_d^{af_0}$ means just ${\cal P}_d^{gr}$ --- the graded counterpart of 
${\cal P}_d$ (of course ${\cal DP}_d^{gr}\simeq {\cal DP}_d$).
The theory of ``$i$--affine functors'' is  parallel but simpler, since $A_i$ is finite dimensional, to
that of ``$\infty$--affine functors''.  In particular, for $j>i$, we have the functors 
\[t^*_{j,i}: {\cal DP}_d^{af_{j}}\ra {\cal DP}_d^{af_i}\]
induced by the embeddings $A_i\subset A_{j}$ and their adjoints
$z^*_{j,i}$. We also consider  the infinite variants:
\[t^*_{\infty,i}: {\cal DP}_d^{af_{\infty}}\ra {\cal DP}_d^{af_i}\]
and their adjoints $z^*_{\infty,i}$.
 We have also the adjunction $\{{\bf C}^{af_i},{\bf K}^{af_i}\}$ between the categories
${\cal DP}_d^{af_i}$ and ${\cal DP}_{dp^i}$.
\begin{lem}
	We have the following isomorphisms of functors:
	\begin{itemize} 
		\item $t^*_{\infty,i}\simeq t^*_{i+1,i}\circ t^*_{\infty,i}$
		\item ${\bf C}^{af_i}\circ t^*_{i+1,i}\simeq {\bf K}\circ {\bf C}^{af_{i+1}}$
		%\item $t^*_{i+1,i}\circ{\bf K}^{af_{i+1}}\simeq{\bf K}^{af_{i}}\circ {\bf K}$.
\end{itemize}
\end{lem}
{\bf Proof:}  The first isomorphism is obvious. In order to get the second one we evaluate both
sides on the cofibrant generator $h^{U\ot A_{i+1}}$ of ${\cal DP}_d^{af_{i+1}}$.
Since $h^{U\ot A_{i+1}}=z^*(\Gamma^{d,U})$, we have
\[{\bf K} ({\bf C}^{af_{i+1}}(h^{U\ot A_{i+1}}))={\bf K} ({\bf K}^{af_{i+1}}(z^*(\Gamma^{d,U})))={\bf K}({\bf C}^{i+1}(\Gamma^{d,U}))=
(\Gamma^{d,U})^{(i+1)}_{A_1},\]
and we recall 
that  \[(\Gamma^{d,U})^{(i+1)}_{A_1}(V):=\Gamma^d(V^{(i)}\ot A_1^{(i+1)}\ot U^*).\]
On the other hand we have
\[{\bf K}^{af_i}(t^*_{i+1,i}(h^{U\ot A_{i+1}}))= {\bf K}^{af_i}(h^{U\ot A_1^{i+1}\ot A_i})=
(\Gamma^{d,U})^{(i+1)}_{A_1}.\]
%The last isomorphism is best seen when we evaluate the both sides on the generator
%$S^{p^{i+1}\la}$ (the both sides are clearly trivial on generators $S^{\mu}$ for 
%$\mu\neq p^{i+1}\la$). Then we have 
%\[t^*_{i+1,i}({\bf K}^{af_{i+1}}(S^{p^{i+1}\la}))=t^*_{i+1,i}(S^{\la})=S^{\la}\]
%where $S^{\la}$ as an object in the affine category means the functor 
%$V\ot A_j\mapsto S^{\la}(V)$ (ie. we kill the $A_j$--structure).
%Similarly we get
%\[{\bf K}^{af_{i}}({\bf K}(S^{p^{i+1}\la}))={\bf K}^{af_{i}}(S^{p^i\la})=S^{\la}.\]
\qed
Now for $F\in{\cal P}_d^{af_{\infty}}$ we consider the collection 
$\widetilde{\gamma}(F)_i:=\{{\bf C}^{af_i}(t^*_{\infty,i}(F))\}$ for $i\geq 0$. 
I claim that $\widetilde{\gamma}(F)$ may be equipped with a structure of ${\bf K}$-spectrum.
Indeed, by using the second and first part of Lemma 6.1 we obtain isomorphisms:
\[
{\bf K}(\widetilde{\gamma}(F)_{i+1})= {\bf K}({\bf C}^{af_{i+1}}(t^*_{\infty,i+1}(F)))\simeq
C^{af_i}(t^*_{i+1,i}(t^*_{\infty,i+1}(F)))\simeq
\]
\[
{\bf C}^{af_i}(t^*_{\infty,i}(F))=\widetilde{\gamma}(F)_i.
\]
Thus $\widetilde{\gamma}$ is a functor from ${\cal P}_d^{af_{\infty}}$
to ${\cal SP}_d$. Then we extend $\widetilde{\gamma}$ to the functor between derived categories:
\[\gamma:{\cal DP}_d^{af_{\infty}}\ra{\cal DSP}_d.\]
by applying $\widetilde{\gamma}$ to  the cofibrant objects in ${\cal KP}_d^{af_{\infty}}$.
\begin{theo} The functor $\gamma$ satisfies the following properties:
\begin{enumerate}
\item There are isomorphisms of functors $\gamma\circ z^*\simeq 
{\bf C}^{\infty}$,  ${\bf f}^{af_{\infty}}\simeq {\bf f}^{st}\circ \gamma$.
\item  $\gamma$ is  a full embedding and it restricts to an equivalence
${\cal DP}_d^{faf_{\infty}}\simeq {\cal DP}_d^{st}$.
\end{enumerate}
\end{theo}
{\bf Remark 6.3 } Theorem 6.2 provides a comparison between the ideas of Sections 
2-3 and Sections 4-5. The first part shows that $\gamma$ is compatible with all our previous constructions. The crucial is the second part, which shows that, in a sense, ${\cal P}_d^{af_{\infty}}$ is a more economical construction than ${\cal SP}_d$, but they become equivalent when restricted to the subcategories consisting of finite objects.\newline
\mbox{}\\
{\bf Proof of Theorem 6.2} First we observe that $t^*_{\infty,i}\circ z^*=(z^*_{i,0})_{A_{\infty}^{(i)}}$. Hence  
for $F\in{\cal KP}_d$ we obtain 
\[\gamma(z^*(F))=\{{\bf C}^{af_i}(t^*_{\infty,i}(z^*(F))) \}=\{{\bf C}^{af_i}(z^*_{i,0}(F_{A_{\infty}^{(i)}}))\}=\{{\bf C}^i(F_{A_{\infty}^{(i)}})\}=F^{(\infty)}_{\Ain}.\] 
Now we recall that by Proposition 4.7 the spectra $F^{(\infty)}_{\Ain}$ and ${\bf C}^{\infty}(F)$ are naturally stably quasi-isomorphic, hence equivalent in 
${\cal DSP}_d$, which shows the first isomorphism.\newline 
In order to establish the second isomorphism we evaluate the both sides on the generator
$h^{U\ot\Ain}$. On the one hand we have \[{\bf f}^{af_{\infty}}(h^{U\ot\Ain})=
{\bf f}(\Gamma^{d,U}).\] On the other hand: \[{\bf f}^{st}(\gamma(h^{U\ot\Ain}))= {\bf f}^{st}((\Gamma^{d,U})^{(\infty)}_{\Ain})\simeq
{\bf f}^{st}({\bf C}^{\infty}(\Gamma^{d,U}))={\bf f}(\Gamma^{d,U}).\]
In order to obtain the second part of Theorem 6.2 we recall that s $h^{U\ot\Ain}=
z^*(\Gamma^{d,U})$, thus  we have 
$\gamma(h^{U\ot\Ain})\simeq \Gamma^{d,U})^{(\infty)}_{\Ain}\simeq {\bf C}^{\infty}(\Gamma^{d,U})$. Hence ,since $ {\bf C}^{\infty}(\Gamma^{d,U})$ form a set of finite  generators of ${\cal DP}_d^{st}$, by [Ke, Lemma 4.2], it suffices to show that $\gamma$ induces bijections:
\[
{\rm Hom}_{{\cal DP}_d^{af_{\infty}}}(h^{U\ot\Ain},h^{W\ot\Ain}[j])\simeq 
{\rm Hom}_{{\cal DSP}_d}((\Gamma^{d,U})^{(\infty)}_{\Ain},(\Gamma^{d,W}[j])^{(\infty)}_{\Ain})
\]
for all spaces $U,W$ and shifts $j$. It is easy to see that the source and target are abstractly isomorphic. However, since $\gamma$ uses $t^*_{\infty,i}$, it may appear that it can kill morphisms. For this reason, let us look carefully how $\gamma$ acts on morphisms.
First of all, since $h^{U\ot\Ain},h^{W\ot\Ain}$ are cofibrant, we have:
\[
{\rm Hom}_{{\cal DP}_d^{af_{\infty}}}(h^{U\ot\Ain},h^{W\ot\Ain}[*])\simeq 
{\rm Hom}_{{\cal P}_d^{af_{\infty}}}(h^{U\ot\Ain},h^{W\ot\Ain})\simeq\] 
\[{\rm Hom}_{{\cal P}_d^{af_{\infty}}}(z^*(\Gamma^{d,U}),z^*(\Gamma^{d,W})).
\]
Now we recall that by Proposition 4.7 $(F)^{\infty}_{\Ain}\simeq \Theta^{\infty}({\bf C}^{\infty}(F))$. Thus by general theory of spectra (and the cofibrance of 
${\bf C}^{\infty}(\Gamma^{d,U})$) we have:
\[
{\rm Hom}_{{\cal DSP}_d}((\Gamma^{d,U})^{(\infty)}_{\Ain},(\Gamma^{d,W}[*])^{(\infty)}_{\Ain})
\simeq
{\rm Hom}_{{\cal DSP}_d}({\bf C}^{\infty}(\Gamma^{d,U}),(\Gamma^{d,W}[*])^{(\infty)}_{\Ain})
 \simeq
 \]
 \[
 {\rm Hom}_{{\cal SP}_d}({\bf C}^{\infty}(\Gamma^{d,U}),(\Gamma^{d,W})^{(\infty)}_{\Ain})
\simeq {\rm Hom}_{{\cal P}_d}(\Gamma^{d,U},(\Gamma^{d,W})_{\Ain})\simeq
{\rm Hom}_{{\cal P}_d}(\Gamma^{d,U},t^*\circ z^*(\Gamma^{d,W})).
\]
We point out that the first  bijection is in general  induced by the canonical map of spectra
${\bf C}^{\infty}(F)\ra \Theta^{\infty}({\bf C}^{\infty}(F))$ which on the zeroth level of spectra in our example may be identified with the unit map $\Gamma^{d,U}\ra t^*\circ z^*(\Gamma^{d,U})$ by Proposition 4.7. Therefore, under these identifications, the action 
of $\gamma$ on the morphisms can be described as the composite:
\[
{\rm Hom}_{{\cal P}_d^{af_{\infty}}}(z^*(\Gamma^{d,U}),z^*(\Gamma^{d,W}))\ra
{\rm Hom}_{{\cal P}_d^{af_{\infty}}}(t^*\circ z^*(\Gamma^{d,U}),t^*\circ z^*(\Gamma^{d,W}))
\ra\]
\[
{\rm Hom}_{{\cal P}_d}(\Gamma^{d,U})),t^*\circ z^*(\Gamma^{d,W}))
\]
where the first map is applying $t^*$ and the second one is induced by the aforementioned unit 
map.  This composite coincides with the adjunction isomorphism by general theory of adjunctions. \qed

\end{document}